\documentclass[preprint,12pt,authoryear]{elsarticle}

\usepackage{amssymb}

\journal
{}

\begin{document}

\begin{frontmatter}

\title{\Large\bf On the number of edges in some graphs \thanksref{support}}
\thanks[support]{Project supported by the National Science Foundation of China (No.61379021; No. 11401290),
 NSF of Fujian (2015J01018; 2016J01027; 2016J01673; 2017J01404; 2018J01423), Fujian Provincial Training Foundation
for "Bai-Quan-Wan Talents Engineering", Project of Fujian Education Department (JZ160455),
the Institute of Meteorological Big Data-Digital Fujian and Fujian Key Laboratory of Data Science and Statistics.}
\author[Lai]{Chunhui Lai\thanksref{correspond}}\ead{laich2011@msn.cn; laichunhui@mnnu.edu.cn}
\address[Lai]{School of Mathematics and Statistics,\\
 Minnan Normal University, Zhangzhou, Fujian,  P.R. China.}

\thanks[correspond]{Corresponding author}

\begin{abstract}

In 1975, P. Erd\H os proposed the problem of determining the
maximum number $f(n)$ of edges in a  graph with $n$ vertices in
which any two cycles are of different
 lengths. The sequence $(c_1,c_2,\cdots,c_n)$  is the cycle length
distribution of a graph $G$ with $n$ vertices,  where  $c_i$ is the
number of cycles of length $i$ in $G$. Let $f(a_1,a_2,\cdots,$
$a_n)$ denote the maximum possible number of edges in a graph which satisfies
$c_i\leq a_i$,  where $a_i$  is a nonnegative integer. In 1991, Shi posed the
problem of determining $f(a_1,a_2,\cdots,a_n)£¬$   which extended
the problem due to Erd\H os. It is clear that $f(n)=f(1,1,\cdots,1)$.
Let $g(n,m)=f(a_1,a_2,\cdots,a_n),$ where $a_i=1$ if $i/m$ is an integer, and $a_i=0$ otherwise.
It is clear that $f(n)=g(n,1)$.
 We prove that
$\liminf\sb {n \to \infty} {f(n)-n \over \sqrt n} \geq \sqrt {2 +
\frac{40}{99}},$   which is better than the previous bounds $\sqrt
2$ (Shi, 1988), and $\sqrt {2 +
\frac{7654}{19071}}$ (Lai, 2017).
We show that $\liminf_{n \rightarrow \infty} {g(n,m)-n\over \sqrt \frac{n}{m}}
  > \sqrt {2.444},$ for all even integers $m$.  We make the following conjecture:
 $\liminf\sb {n \to \infty} {f(n)-n \over \sqrt n} > \sqrt {2.444}.$
\par
 \par

 \par
\end{abstract}

\begin{keyword}
Graph,  cycle,  number of edges.
\def\MSC{\par\leavevmode\hbox{\em AMS 2000 MSC:\ }}%
\MSC 05C38, 05C35.
\end{keyword}
\end{frontmatter}


\section{Introduction}
\label{1}
\par
Let $f(n)$ be the maximum number of edges in a graph with $n$ vertices
in which no two cycles have the same length. In 1975, Erd\H os
raised the problem of determining $f(n)$ (see Bondy and
Murty~\cite{J.A. Bondy and U.S.R. Murty (1976)}, p.247, Problem 11).
Shi~\cite{Y. Shi (1988)} proved a lower bound.
\par
\par
\noindent{\bf Theorem 1 (Shi~\cite{Y. Shi (1988)})}
$$f(n)\geq n+
 [(\sqrt {8n-23} +1)/2]$$ for $n\geq 3$.
\par
 Chen, Lehel, Jacobson and
 Shreve~\cite{G. Chen J. Lehel M. S. Jacobson and W. E. Shreve (1998)}, Jia~\cite{X. Jia (1996)},
 Lai~\cite{C. Lai (1993),C. Lai (2001), C. Lai (2003)},
 Shi~\cite{Y. Shi (1992),Y. Shi (1994)} obtained some additional related results.
\par
 Boros, Caro, F\"uredi and Yuster~\cite{E. Boros Y. Caro Z. Furedi and R. Yuster (2001)}
 proved an upper bound as follows.
\par
\noindent{\bf Theorem 2 (Boros, Caro, F\"uredi and Yuster~\cite{E.
Boros Y. Caro Z. Furedi and R. Yuster (2001)})}
  For $n$ sufficiently large, $$f(n) < n+1.98\sqrt{n}.$$
\par
Lai~\cite{C. Lai (2017)} improved the lower bound by Shi as follows.
\par
\par
\noindent{\bf Theorem 3 (Lai~\cite{C. Lai (2017)})}
  Let $ t=1260r+169 \,\ (r\geq 1)$, then $$f(n)\geq n+\frac{107}{3}t+\frac{7}{3}$$ for
  $n\geq
  \frac{2119}{4}t^{2}+87978t+\frac{15957}{4}$.
  \par
 \par
 Lai~\cite{C. Lai (1993)} proposed the following conjecture:

 \par
\noindent{\bf Conjecture 4 (Lai~\cite{C. Lai (1993)})}
$$\liminf_{n \rightarrow \infty} {f(n)-n\over \sqrt n} \leq \sqrt {3}.$$
\par

It would be nice to prove that
 \par
 $$\liminf_{n \rightarrow \infty} {f(n)-n\over \sqrt n} \leq \sqrt {3+ \frac{3}{5}}.$$
\par
 Survey papers on this problem can be found in Tian~\cite{F. Tian (1986)},
 Zhang~\cite{K. Zhang (2007)}, Lai and Liu~\cite{C. Lai (2014)}.
\par
The progress of all 50 problems in~\cite{J.A. Bondy and U.S.R. Murty
(1976)} can be found in Locke~\cite{S. C. Locke (2016)}.
\par
\par
The sequence $(c_1,c_2,\cdots,c_n)$  is the cycle length
distribution of a graph $G$ with $n$ vertices,  where  $c_i$ is the
number of cycles of length $i$ in $G$. Let $f(a_1,a_2,\cdots,$
$a_n)$ denote the maximum possible number of edges in a graph which satisfies
$c_i\leq a_i$,  where $a_i$  is a nonnegative integer. Shi~\cite{Y. Shi (1991)} posed the
problem of determining $f(a_1,a_2,\cdots,a_n)£¬$   which extended
the problem due to Erd\H os. It is clear that $f(n)=f(1,1,\cdots,1)$.
Let $g(n,m)=f(a_1,a_2,\cdots,a_n),$ where $a_i=1$ if $i/m$ is an integer, and $a_i=0$ otherwise.
It is clear that $f(n)=g(n,1)$.
 \par

 In this paper, we obtain the following results.

   \par
\bigskip
  \noindent{\bf Theorem 5} Let $m$ be even, $s_{1}> s_{2}, s_{1}+3s_{2}>k$ , then $$g(n,m)\geq n+(k+s_{1}+2s_{2}+1)t-1$$ for
 $n \geq (\frac{3}{4}mk^{2}+\frac{1}{2}mks_{1}+\frac{3}{2}mks_{2}+\frac{1}{2}ms_{1}^{2}+\frac{3}{2}ms_{1}s_{2}+\frac{9}{4}ms_{2}^{2}+mk+ms_{1}+3ms_{2}+\frac{1}{2}m)t^{2}+
 (\frac{1}{4}mk+\frac{1}{2}ms_{1}+\frac{3}{4}ms_{2}-k-s_{1}-2s_{2}+\frac{1}{2}m-1)t+1$.
 \par
  \par
\bigskip
  \noindent{\bf Theorem 6} Let $ t=1260r+169 \,\ (r\geq 1)$, then $$f(n)\geq n+\frac{119}{3}t-\frac{26399}{3}$$for
 $n \geq \frac{1309}{2}t^{2}-\frac{1349159}{6}t+\frac{6932215}{3}$.

 \par
  \section{Proof of Theorem 5}
\label{2}
\par
 {\bf Proof.} Let
  $n_{t}=(\frac{3}{4}mk^{2}+\frac{1}{2}mks_{1}+\frac{3}{2}mks_{2}+\frac{1}{2}ms_{1}^{2}+\frac{3}{2}ms_{1}s_{2}+\frac{9}{4}ms_{2}^{2}+mk+ms_{1}+3ms_{2}+\frac{1}{2}m)t^{2}+
 (\frac{1}{4}mk+\frac{1}{2}ms_{1}+\frac{3}{4}ms_{2}-k-s_{1}-2s_{2}+\frac{1}{2}m-1)t+1$, $m$ be even, $s_{1}> s_{2}, s_{1}+3s_{2}>k$, $n\geq n_{t}.$
  It suffice to show that there exists a graph $G$ on $n$ vertices with $ n+(k+s_{1}+2s_{2}+1)t-1$ edges such
  that all cycles in $G$ have distinct lengths and all the lengths of cycles are the multiple of $m$.

  Now we construct the graph $G$ which consists of a number of subgraphs: $B_i$,
  ($0\leq i\leq s_{1}t,i=s_{1}t+j$ $(1\leq j \leq s_{2}t),$ $i=s_{1}t+s_{2}t+j$ $(1\leq j \leq t)$).

  Now we define these $B_is$. These subgraphs all only have a common vertex $x$,
  otherwise their vertex sets are pairwise disjoint.
  \par
For $1 \leq i\leq s_{2}t,$ let the subgraph $B_{s_{1}t+i}$
consists of a cycle $$xa_i^1a_i^2...a_i^{ms_{1}t+2ms_{2}t+mi-1}x$$  and a path:

  $$xa_{i,1}^1a_{i,1}^2...a_{i,1}^{\frac{ms_{1}t-ms_{2}t+mi}{2}-1}a_i^{\frac{ms_{1}t+ms_{2}t+mi}{2}}.$$
  \par
  Based on the construction,  $B_{s_{1}t+i}$ contains exactly three
cycles of lengths:
$$ms_{1}t+mi, ms_{1}t+ms_{2}t+mi, ms_{1}t+2ms_{2}t+mi.$$
\par

   For $1\leq i\leq t,$ let the subgraph  $B_{s_{1}t+s_{2}t+i}$ consists of a
  cycle $$C_{s_{1}t+s_{2}t+i}=xy_i^1y_i^2...y_i^{ms_{1}t+3ms_{2}t+mk(k+1)t+mi-1}x$$  and $k$ paths
  sharing a common vertex $x$, the other end vertices are on the
  cycle $C_{s_{1}t+s_{2}t+i}$:

  $$xy_{i,p}^1y_{i,p}^2...y_{i,p}^{\frac{ms_{1}t+3ms_{2}t-mkt+m(p-1)t+mi}{2}-1}y_i^{\frac{ms_{1}t+3ms_{2}t+mk(2p-1)t+m(p-1)t+mi}{2}} (p=1,2,...,k).$$

\par
As a cycle with $k$ chords contains ${{k+2} \choose 2}$ distinct
cycles, $B_{s_{1}t+s_{2}t+i}$ contains exactly $\frac{(k+2)(k+1)}{2}$ cycles of lengths:
 \par
 $$ms_{1}t+3ms_{2}t+mkht+(h+j-1)mt+mi (j\geq 1, h\geq 0, k+1\geq j+h).$$
 \par
 $B_{0}$ is a path with an end vertex $x$ and length $n-n_{t}$. The other $B_i$ is
 simply a cycle of length $mi$.
  \par Then  $g(n,m)\geq n+(k+s_{1}+2s_{2}+1)t-1,$ for $n\geq n_{t}.$
 \par

This completes the proof.
  \par
   From Theorem 5, we have $$\liminf_{n \rightarrow \infty} {g(n,m)-n\over \sqrt \frac{n}{m}}
  \geq$$ $$ \sqrt {\frac{(k+s_{1}+2s_{2}+1)^{2}}{(\frac{3}{4}k^{2}+\frac{1}{2}ks_{1}+\frac{3}{2}ks_{2}+\frac{1}{2}s_{1}^{2}+\frac{3}{2}s_{1}s_{2}+\frac{9}{4}s_{2}^{2}+k+s_{1}+3s_{2}+\frac{1}{2})}},$$
for all even integers $m$.
\par
  Let $s_{1}=28499066, s_{2}=4749839, k=14249542$, then
  $$\liminf_{n \rightarrow \infty} {g(n,m)-n\over \sqrt \frac{n}{m}}
  > \sqrt {2.444},$$ for all even integers $m$.
  \vskip 0.2in

 \par
  \section{Proof of Theorem 6}
\label{3}
\par
 {\bf Proof.} Let
  $ n_{t}=\frac{1309}{2}t^{2}-\frac{1349159}{6}t+\frac{6932215}{3},$ $t=1260r+169,r\geq 1,$ $n\geq n_{t}.$
  It suffice to show that there exists a graph $G$ on $n$ vertices with $ n+\frac{119}{3}t-\frac{26399}{3}$ edges such
  that all cycles in $G$ have distinct lengths.

  Now we construct the graph $G$ which consists of a number of subgraphs: $B_i$,
  ($0\leq i\leq 22t,
  i=22t+j$ $(1\leq j \leq \frac{5t-8}{3}),$
  $i=23t+\frac{2t-2}{3}+j$ $(1\leq j \leq \frac{5t-8}{3}),$
  $i=32t+j-60$ $(58\leq j \leq t-742)$).

  Now we define these $B_is$. These subgraphs all only have a common vertex $x$,
  otherwise their vertex sets are pairwise disjoint.
  \par
For $1 \leq i\leq \frac{5t-8}{3},$ let the subgraph $B_{22t+i}$
consists of a cycle $$xa_i^1a_i^2...a_i^{28t+\frac{2t-2}{3}+2i-3}x$$  and a path:

  $$xa_{i,1}^1a_{i,1}^2...a_{i,1}^{\frac{56t-2}{6}}a_i^{\frac{76t-4}{6}+i}.$$
  \par
  Based on the construction,  $B_{22t+i}$ contains exactly three
cycles of lengths:
$$22t+i, 25t+\frac{t-1}{3}+i-1, 28t+\frac{2t-2}{3}+2i-2.$$
\par
   \par
For $1 \leq i\leq \frac{5t-8}{3},$ let the subgraph $B_{23t+\frac{2t-2}{3}+i}$
consists of a cycle $$xb_i^1b_i^2...b_i^{28t+\frac{2t-2}{3}+2i-2}x$$  and a path:

  $$xb_{i,1}^1b_{i,1}^2...b_{i,1}^{11t-1}b_i^{\frac{76t-4}{6}+i}.$$
  \par
  Based on the construction,  $B_{23t+\frac{2t-2}{3}+i}$ contains exactly three
cycles of lengths:
$$23t+\frac{2t-2}{3}+i, 27t+i-1, 28t+\frac{2t-2}{3}+2i-1.$$
\par
   For $58\leq i\leq t-742,$ let the subgraph  $B_{32t+i-60}$ consists of a
  cycle $$C_{32t+i-60}=xy_i^1y_i^2...y_i^{137t+11i+890}x$$  and ten paths
  sharing a common vertex $x$, the other end vertices are on the
  cycle $C_{32t+i-60}$:

  $$xy_{i,1}^1y_{i,1}^2...y_{i,1}^{11t-2}y_i^{21t-59+i}$$
  $$xy_{i,2}^1y_{i,2}^2...y_{i,2}^{12t-2}y_i^{31t-53+2i}$$
  $$xy_{i,3}^1y_{i,3}^2...y_{i,3}^{12t-2}y_i^{41t+156+3i}$$
  $$xy_{i,4}^1y_{i,4}^2...y_{i,4}^{13t-2}y_i^{51t+155+4i}$$
  $$xy_{i,5}^1y_{i,5}^2...y_{i,5}^{13t-2}y_i^{61t+155+5i}$$
  $$xy_{i,6}^1y_{i,6}^2...y_{i,6}^{14t-2}y_i^{71t+154+6i}$$
  $$xy_{i,7}^1y_{i,7}^2...y_{i,7}^{14t-2}y_i^{81t+153+7i}$$
  $$xy_{i,8}^1y_{i,8}^2...y_{i,8}^{15t-2}y_i^{91t+147+8i}$$
  $$xy_{i,9}^1y_{i,9}^2...y_{i,9}^{15t-2}y_i^{101t+149+9i}$$
  $$xy_{i,10}^1y_{i,10}^2...y_{i,10}^{16t-2}y_i^{111t+151+10i}.$$

\par
As a cycle with $d$ chords contains ${{d+2} \choose 2}$ distinct
cycles, $B_{32t+i-60}$ contains exactly 66 cycles of lengths:
 \par
 $$\begin{array}{llll}32t+i-60,& 33t+i+4,& 34t+i+207,& 35t+i-3,\\
  36t+i-2,& 37t+i-3,& 38t+i-3,& 39t+i-8,\\
  40t+i,& 41t+i,& 42t+i+739,& 43t+2i-54,\\
  43t+2i+213,& 45t+2i+206,& 45t+2i-3,& 47t+2i-3,\\
  47t+2i-4,& 49t+2i-9,& 49t+2i-6,& 51t+2i+2,\\
  51t+2i+741,& 53t+3i+155,& 54t+3i+212,& 55t+3i+206,\\
  56t+3i-4,& 57t+3i-4,& 58t+3i-10,& 59t+3i-7,\\
  60t+3i-4,& 61t+3i+743,& 64t+4i+154,& 64t+4i+212,\\
  66t+4i+205,& 66t+4i-5,& 68t+4i-10,& 68t+4i-8,\\
  70t+4i-5,& 70t+4i+737,& 74t+5i+154,& 75t+5i+211,\\
  76t+5i+204,& 77t+5i-11,& 78t+5i-8,& 79t+5i-6,\\
  80t+5i+736,& 85t+6i+153,& 85t+6i+210,& 87t+6i+198,\\
  87t+6i-9,& 89t+6i-6,& 89t+6i+735,& 95t+7i+152,\\
  96t+7i+204,& 97t+7i+200,& 98t+7i-7,& 99t+7i+735,\\
  106t+8i+146,& 106t+8i+206,& 108t+8i+202,& 108t+8i+734,\\
  116t+9i+148,& 117t+9i+208,& 118t+9i+943,& 127t+10i+150,\\
  127t+10i+949,& 137t+11i+891.&&\end{array}$$
 \par
 $B_{0}$ is a path with an end vertex $x$ and length $n-n_{t}$. The other $B_i$ is
 simply a cycle of length $i$.
  \par Then  $f(n)\geq n+\frac{119}{3}t-\frac{26399}{3},$ for $n\geq n_{t}.$
 \par

This completes the proof.
  \par
  \vskip 0.2in

  From Theorem 6, we have $$\liminf_{n \rightarrow \infty} {f(n)-n\over \sqrt n}
  \geq \sqrt {2 +
\frac{40}{99}},$$
  which is better than the previous bounds $\sqrt 2$ (see ~\cite{Y. Shi (1988)}), and $\sqrt {2 +
\frac{7654}{19071}}$
  (see ~\cite{C. Lai (2017)}).

 \vskip 0.2in
  \par
  Combining this with Boros, Caro, F\"uredi and Yuster's upper bound, namely Theorem 2, we get
  $$1.98\geq \limsup_{n \rightarrow \infty} {f(n)-n\over \sqrt n} \geq
  \liminf_{n \rightarrow \infty} {f(n)-n\over \sqrt n}\geq \sqrt {2 +
\frac{40}{99}}.$$
  \par
  From the proof of Theorem 6, we have
 $$\liminf_{n \rightarrow \infty} {g(n,m)-n\over \sqrt \frac{n}{m}}
  \geq \sqrt {2 + \frac{40}{99}},$$
  for all integers $m$.
  \par
 If $m=1,$ $1\leq i\leq t,$ there exists the subgraph similar to  $B_{s_{1}t+s_{2}t+i}$ consists of a
  cycle $C_{s_{1}t+s_{2}t+i}$ and $k$ paths
  sharing a common vertex $x$, the other end vertices are on the
  cycle $C_{s_{1}t+s_{2}t+i}$ such
  that all cycles in $B_{s_{1}t+s_{2}t+i}$ have distinct lengths,
  then we could obtain $$\liminf\sb {n \to \infty} {f(n)-n \over \sqrt n} > \sqrt {2.444}> \sqrt {2 +
\frac{40}{99}}.$$
  But we only for $m=1,$ $58\leq i\leq t-742,$ construct a subgraph similar to  $B_{s_{1}t+s_{2}t+i}$ consists of a
  cycle $C_{s_{1}t+s_{2}t+i}$ and ten paths
  sharing a common vertex $x$, the other end vertices are on the
  cycle $C_{s_{1}t+s_{2}t+i}$ such
  that all cycles in $B_{s_{1}t+s_{2}t+i}$ have distinct lengths
  and obtain $$\liminf\sb {n \to \infty} {f(n)-n \over \sqrt n} \geq \sqrt {2 +
\frac{40}{99}}.$$
\par
 \par
 We make the following conjecture:
  \par
  \bigskip
  \noindent{\bf Conjecture 7} $$\liminf\sb {n \to \infty} {f(n)-n \over \sqrt n} > \sqrt {2.444}.$$
\par

\par
  \section{Acknowledgment}
\label{4}
\par
The author would like to thank Professor Endre Boros, Yair Caro, Ronald J. Gould, Gyula O.H. Katona for their advice.
The author would like to thank the referees for their many valuable comments and suggestions.
\par





\end{document}